\documentclass[11pt,a4paper]{article}

\usepackage[T2A]{fontenc}
\usepackage[utf8]{inputenc}
\usepackage[english]{babel}

\usepackage{amssymb,amsthm,amsmath,amsfonts,array,setspace}
\usepackage[text={15cm,25cm}, includehead]{geometry}

\theoremstyle{plain}
\newtheorem{theorem}{Theorem}         %для наших теорем
\newtheorem*{primer}{Example}

\begin{document}
\pagestyle{myheadings}

\makeatletter
\renewcommand{\@oddhead}{\hfil 
    A set of 12 numbers is not determined by its set of 4-sums
\hfil\thepage}
\makeatother

\title{A set of 12 numbers is not determined by its set of 4-sums}
\author{J. E. Isomurodov\thanks{ITMO University, St.Petersburg, Russia.} 
\and 
K. P. Kokhas\thanks{St.Petersburg State University, St.Petersburg, Russia} }

\date{}

\maketitle

\begin{abstract}
We present two sets of 12 integers that have the same sets of 4-sums. The proof of the fact that a
set of 12 numbers is uniquely determined by the set of its 4-sums published 50 years ago is wrong,
and we demonstrate an incorrect calculation in it
\end{abstract}

\subsection{Introduction}
For an arbitrary set $A$ of $n$ real numbers consider a set $A^{(k)}$ of all sums of
$k$ distinct elements of $A$. We call this set the \emph{set of $k$-sums}. 
For which pairs $(n, k)$ the initial set $A$ is uniquely determined by the set of its
$k$-sums?  This question was posed by L.\,Moser \cite{Moser1957} in American Math.\ Monthly in 1957. 
It turned out to be highly nontrivial and  was discussed in several papers in 1960-th \cite{gordon}, \cite{straus},
\cite{ewell} and 1990-th \cite{fomin}, \cite{boman}. In the book 
``Unsolved problems in number theory'' by {\it R.\,Guy} this
problem is mentioned as one of the unsolved problems in additive number theory \cite[problem C5]{Guy1994}.

In 1958 J.\,L.\,Selfridge and E.\,G.\,Straus \cite{straus} proved that in the
case $(n, 2)$ the initial set $A$ can be recovered from $A^{(k)}$ if and only if $n$ is not a power of 2.
For the general case, they constructed a system of Diophantine equations on the variables
$n$ and $k$. If a pair $(n, k)$ satisfies neither of these equations (this is a typical case), 
the initial set $A$ can be recovered  uniquely. 
But if a pair $(n, k)$ satisfies at least one of the equations
(we have exactly this case for $(n,k)=(12, 4)$), the possibility of recovery is not clear. 
In 1962 A.\,S.\,Fraenkel, B. Gordon, E.\,G.\,Straus
\cite{gordon} showed that it is sufficient to consider the problem over integers
and that for every $k$ there exists only a finite set of numbers $n$ for which 
$A$ is not uniquely determined by $A^{(k)}$. 
For <<bad>> cases $(8,2)$, $(16,2)$, $(6,3)$ and
$(12,4)$, they gave an upper bound on the number of sets $A$ for which the sets $A^{(k)}$ coincide. 
In 1968 J.\,Ewell \cite{ewell} proved that for the case $(6, 3)$ 
the maximal number of sets $A$ with the same set of $k$-sums is 4.
In 1994 D.\,V.\,Fomin and O.\,T.\,Izhboldin \cite{fomin} 
showed, for cases $(27, 3)$ and $(486, 3)$, the examples of two different sets with
the same set of $k$-sums. 
Similar examples were constructed independently by J.\,Boman and S.\,Linusson \cite{boman}.

J.\,Ewell \cite{ewell} proved also that for the case $(12, 4)$
the initial set $A$ can be recovered uniquely.
Unfortunately, this is not correct, and below we show a mistake in Ewell's proof. 
While we were trying to fix the mistake using the computer
algebra system Maple, we found an example of two sets of 12 numbers with
the same set of 4-sums. So, the following theorem holds.

\begin{theorem}
There exist two sets of 12 integers that have the same set of 4-sums.
\end{theorem}

There do not exist three sets of 12 numbers with the same set of 4-sums. This statement is
proven in~\cite{gordon}. It follows from our calculations that
the example of two sets  with the same set of 4-sums is \emph{unique} (up to shifts and scalings).

\subsection{The original proof}

Let $A$ be an arbitrary set of $n$ (not necessarily different) numbers 
$$
A = \{a_1, a_2, \ldots, a_n\}. 
$$ 
A \emph{$k$-sum} is a number of the form
\begin{equation*}
    a_{i_1} + a_{i_2} + \ldots + a_{i_k}
    \text{, where } 1 \leqslant i_1 < i_2 < \ldots < i_k \leqslant n.
\end{equation*}
Denote by $A^{(k)}$ the set of all $k$-sums. We consider
the special case $n=12$, $k=4$ of the following main problem on the sets~$A^{(k)}$.

\medskip 
\centerline{ \it Find all pairs $(n, k)$ such that the set $A$ is determined by
    the set $A^{(k)}$} 
\medskip

We briefly describe the approach to this problem developed in \cite{straus, gordon, ewell}.

In \cite{gordon}, it is shown that the answer to the question does not depend on 
whether the numbers $a_i$ are integer, real, or complex. 
So we assume that our sets contain complex numbers. 
Consider the symmetric functions of set $A$. 
For $1 \leq j \leq n$  and for any set of nonnegative integers $\{p_1, p_2, \ldots, p_j\}$ 
denote by $S_{p_1, p_2, \ldots, p_j}$ the symmetric monomial power sum 
\begin{equation} \label{eqn:monomialpowersums}
S_{p_1, p_2, \ldots, p_j} = \sum_{D(j)} a_{i_1}^{p_1} a_{i_2}^{p_2} \ldots
a_{i_j}^{p_j},
\end{equation}
where the summation runs over the set $D(j)$ of all ordered $j$-subsets  $\{i_1, i_2 , \ldots,i_j\}
\subset \{1, \ldots, n\}$. In particular \[ S_p = \sum_{i = 1}^{n} a_i^p. \]

The following observation allows us to reduce the complexity of calculations. It is easy
to check that the shift $\{a_1,\dots, a_n\}\mapsto \{a_1+t,\dots, a_n+t\}$ and the
scaling $\{a_1,\dots, a_n\}\mapsto \{ta_1,\dots, ta_n\}$ do not affect the
solvability of the original problem. Hence, we can put $S_1=0$, $S_2=1$ in the most cumbersome calculations.

It is well known that the set of power sums $S_p$, $p \in \{1, 2 \dots, n\}$ determines the
set $A$ uniquely.  The sum $S_m$ for $m > n$ can be expressed via $S_1$, \dots, $S_n$ 
by repeatedly applying the following formula (see MacMahon's book \cite[page 6]{MacMahon}):
\begin{equation} \label{eqn:McMahon-reduction}
\frac{1}{m} S_m = \sum_{\substack{1p_1 + 2p_2 + \ldots = m \\ p_m = 0}}
(-1)^{^{\sum p_i}} \frac{S_1^{p_1} S_2^{p_2} \ldots}{1^{p_1} 2^{p_2} \ldots
p_1! p_2! \ldots}.
\end{equation}
The monomial power sum $S_{p_1, p_2, \ldots, p_j}$ can also be expressed via $S_1$, \dots,
$S_n$ by applying the following formula:
\begin{equation}
    \label{eqn:sdef} S_{p_1, p_2, \ldots, p_j} = S_{p_1, p_2, \ldots,
    p_{j-1}}S_{p_j} - S_{p_1 + p_j, p_2, \ldots, p_{j - 1}} - \ldots - S_{p_1,
    p_2, \ldots, p_{j - 1} + p_j},
\end{equation}
which follows from the definition of the monomial power sum \eqref{eqn:monomialpowersums}.

Now consider the symmetric functions of set $A^{(k)}$. Let $E_p$ denote the $p$-th
power sum of $A^{(k)}$. Expand the brackets in the definition of $E_p$ (for the case $k=4$):
\begin{equation}\label{eqn:bigwsmall}
    4! E_p = \sum_{D(4)}(a_{i_1} + a_{i_2} + a_{i_3} + a_{i_4})^p
    = \sum_{ p_1 + p_2 + p_3 + p_4 = p}
        \frac{p!}{p_1! p_2! p_3! p_4!} S_{p_1, p_2, p_3, p_4}.
\end{equation}
We can apply formulas \eqref{eqn:sdef}, \eqref{eqn:McMahon-reduction} to the
right-hand side of the last formula and obtain the following expressions 
of $E_p$ via $S_1$, \dots, $S_n$ 
(these expression are given in \cite{ewell}).

{%\footnotesize
\begin{align} 
    E_{1} &= 0 \\
    E_{2} &= 120S_{2} \label{eqn:key2} \\
    E_{3} &= 48S_{3} \label{eqn:key3}\\
    E_{4} &= -48S_{4} +84S_{2}^2 	 \allowdisplaybreaks	\\
    E_{5} &= -120S_{5} +140S_{2}S_{3} \label{eqn:key5}\\ 
    E_{6} &= 0S_{6}+40S_{3}^2-120S_{2}S_{4} +90S_{2}^3 \label{eqn:key6} \\	 
    E_{7} &= 648S_{7}-714S_{2}S_{5}-350S_{3}S_{4} +420S_{2}^2S_{3} \\
    E_{8} &= 1632S_{8}-896S_{2}S_{6}-1120S_{3}S_{5}-280S_{4}^2 
        +560S_{2}S_{3}^2 +105S_{2}^4  \\
    E_{9} &= -3480S_{9} +4824S_{2}S_{7} +1176S_{3}S_{6} +1764S_{4}S_{5} 
        -3024S_{2}^2S_{5} \notag \\
        &\qquad\qquad -2520S_{2}S_{3}S_{4} +1260S_{2}^3S_{3} \allowdisplaybreaks \\
    E_{10} &= -59520S_{10} +42840S_{2}S_{8} +29280S_{3}S_{7} +23520S_{4}S_{6}
        -15120S_{2}^2S_{6} -8400S_{3}^2S_{4}\notag\\
        &\qquad\qquad +3150S_{2}^3S_{4}-9450S_{2}S_{4}^2
        +6300S_{2}^2S_{3}^2 -25200S_{2}S_{3}S_{5}+12600S_{5}^2 \allowdisplaybreaks \\
    E_{11} &= -407352S_{11} +222530S_{2}S_{9} +196350S_{3}S_{8} 
        +155100S_{4}S_{7}-120120S_{2}S_{3}S_{6}\notag\\
        &\qquad\qquad +150612S_{5}S_{6}
        -97020S_{2}S_{4}S_{5}
        -55440S_{3}^2S_{5}+6930S_{2}^3S_{5}\notag\\
        &\qquad\qquad -55440S_{2}^2S_{7}-46200S_{3}S_{4}^2
        +34650S_{2}^2S_{3}S_{4}+15400S_{2}S_{3}^3 \\
    E_{12} &= -2203488S_{12} +964128S_{2}S_{10} +998800S_{3}S_{9}
        +827640S_{4}S_{8}-178200S_{2}^2S_{8} \notag\\
        &\qquad\qquad-459360S_{2}S_{3}S_{7}+744480S_{5}S_{7}
        +373296S_{6}^2-415800S_{2}S_{4}S_{6} \notag\\
        &\qquad\qquad-258720S_{3}S_{3}S_{6}+13860S_{2}^3S_{6}
        -182952S_{2}S_{5}^2-443520S_{3}S_{4}S_{5} \notag\\
        &\qquad\qquad+83160S_{2}^2S_{3}S_{5} -69300S_{4}^3
        +51975S_{2}^2S_{4}^2 +138600S_{2}S_{3}^2S_{4} +15400S_{3}^4
        \label{eqn:key12}
\end{align}
}

Observe that the coefficient of $S_6$ in \eqref{eqn:key6} vanishes. 
This explains why the case (12, 4) is ``bad''. 
If the coefficient were not equal to 0, we could successively find $S_2$,  \dots, $S_{12}$ and 
recover the initial set $A$ uniquely.

Since the equation \eqref{eqn:key6} turned out to be useless, J.\,Ewell added one more
equation to the system:
\begin{align}
E_{14} &= -48517440S_{14} +14260792S_{2}S_{12} +18521776S_{3}S_{11}
        +17649632S_{4}S_{10}\notag\\
        &\qquad\qquad -1513512S_{2}^2S_{10}-5005000S_{2}S_{3}S_{9}
        +15095080S_{5}S_{9}+14030016S_{6}S_{8}\notag\\
        &\qquad\qquad-5675670S_{2}S_{4}S_{8}-3723720S_{3}^2S_{8}
        +45045S_{2}^3S_{8} +7008144S_{7}^2\notag\\
        &\qquad\qquad-5045040S_{2}S_{5}S_{7}-7687680S_{3}S_{4}S_{7}
        -2270268S_{2}S_{6}^2+360360S_{2}^2S_{3}S_{7}\notag\\
        &\qquad\qquad-6726720S_{3}S_{5}S_{6}-3783780S_{4}^2S_{6}
        +630630S_{2}^2S_{4}S_{6} +840840S_{2}S_{3}^2S_{6}\notag\\
        &\qquad\qquad-3531528S_{4}S_{5}^2+378378S_{2}^2S_{5}^2
        +2522520S_{2}S_{3}S_{4}S_{5}+560560S_{3}^3S_{5} \notag \\
        &\qquad\qquad+525525S_{2}S_{4}^3 +1051050S_{3}^2S_{4}^2.
        \label{eqn:key14}
\end{align}
Express the sum $S_{14}$ in terms $S_2$,  \dots, $S_{12}$ by formula
\eqref{eqn:McMahon-reduction} and substitute this expression into the above equation.
We obtain a system where the number of the variables $S_i$ equals to the number of equations. 
Now we will successively eliminate the variables $S_2$, \dots, $S_5$,
$ S_7$, \dots, $S_{12}$ from the system.
To do this, we express the variables $S_2$, \dots, $S_5$ in terms of $E_i$ 
directly from  \eqref{eqn:key2}--\eqref{eqn:key5}, 
and the variables  $S_7$, \dots, $S_{12}$ are
expressed as polynomials in $S_6$ with coefficients depending on $E_i$. 
Thus we obtain the following equation:
\begin{multline} 
    E_{14} = \bigg( \frac{73458}{5465}E_2 \bigg)S_6^2 +
    \bigg(
        \frac{22556178701}{5315943600}E_3 E_5-\frac{889}{12}E_8
        -\frac{15211}{13392}E_4^2 \\
        +\underline{\underline{\frac{4783550233}{119441640960}}}E_2^2 E_4
        -\underline{\frac{9881683541849}{418343497545600}}E_2 E_3^2
        -\underline{\frac{72629302403}{477766563840000}}E_2^4
        \bigg)S_6 + \ldots. \label{eqn:quadroE}
\end{multline}

The last term, denoted by the dots, depends on the parameters $E_i$ only. 
The underlined coefficients presented here are the results of our calculations,
they do not coincide with the coefficients from the paper by J.\,Ewell \cite{ewell}. 
The exact values of the coefficients are not really important,
but J.~Ewell relies on the fact that the doubly underlined coefficient is sufficiently large, 
which ensures that the equation has only one positive root. 
But in our calculations, this coefficient is approximately 3000 times less than in J.\,Ewell's paper. 
The following simple example shows that that the equation can have two positive roots.

\begin{primer} \label{exmp:kp}
    Let $A = \{-1, 0^{10}, 1\}$, then $A^{(4)} = \{-1^{120}, 0^{255}, 1^{120}
    \}$. Then $E_i = 0$ for odd $i$ and $E_i = 240$ for even $i$. Substituting these numbers into
    \eqref{eqn:quadroE}, we obtain an equation with two positive roots, $S_6 = 2$ and $S_6
    = 377762/44361$, which contradicts \cite[theorem 2]{ewell}.
\end{primer}

\subsection{Proof of theorem 1}

Consider the sets
$$
A'=\{0, 0, 1, -1,   2, -2, 4, -4, 7, -7, 7, -7\}
$$
and
$$
A''=\{1, -1, 2, -2, 3, -3,  4, -4, 5, -5,  8, -8\}.
$$
An easy computation shows that the sets $(A')^{(4)}$ and $(A'')^{(4)}$
coincide.

\subsection{Some details}

In this section, we explain how we found the above two sets. Computations were
performed in the computer algebra system Maple.

We try to find the sets $A'=\{a_1',\dots, a_{12}'\}$ and
$A''=\{a_1'',\dots, a_{12}''\}$ whose sets of 4-sums coincide. 
Let $S_p'$, $S_p''$ be the $p$-th powers sums of the sets $A'$ and $A''$, 
and $E_p$ be the $p$-th power sum of the set 
$(A')^{(4)}=(A'')^{(4)}$. Using the equations \eqref{eqn:key2}, \eqref{eqn:key3}, \dots, we express $E_2$,
$E_3$, \dots\ in terms of $S_1'$, \dots, $S_{12}'$.
According to the above remark, we may assume that $S_1'=0$, $S_2'=1$. 

Consider the system of equations \eqref{eqn:key2}, \eqref{eqn:key3}, \dots once again
with respect to new variables $S_i$.
Let us number the equations by their left-hand sides:
the \textit{$k$-th equation} is the equation whose left-hand side contains $E_k$. 
So \eqref{eqn:key2} is the second equation and
\eqref{eqn:key3} is the 12th equation etc. 
In what follows, we use the equations with numbers from 1 to 26.
In them, we substitute the expressions for $E_i$ obtained in the previous paragraph.

Thus we have 26 equations with respect to the variables $S_i$, 
and the left-hand sides $E_i$ of these equations are expressed in terms $S_1'$, \dots, $S_{12}'$. 
Now apply  the same transformations as in the paper by J.\,Ewell for the equations from 13 to 26: 
express $S_{13}$, $S_{14}$, \dots,  $S_{26}$ in terms of $S_2$,  \dots, $S_{12}$ using 
formula~\eqref{eqn:McMahon-reduction} and then  eliminate $S_2$, \dots, $S_5$, $S_7$, \dots,
$S_{12}$ using equations \eqref{eqn:key2}--\eqref{eqn:key12}. 
For example, 14th equation becomes a quadratic equation in $S_6$, 
the only difference with \eqref{eqn:quadroE} is that all $E_i$ here are expressed via $S_j'$.  
The roots of the 14th equation should be equal to $S_6'$ and $S_6''$. 
Indeed, the equation vanishes when we substitute $S_6=S_6'$, and due to this fact 
the second possible root $S_6=S_6''$ can be expressed via $S_1'$, \dots,
$S_{12}'$ (here we put $S_1=0$, but the substitution $S_2=1$ is not yet performed):
\begin{multline*}
    S_6''= -\frac{556877605}{796368672}
    {S_2'}^3+\frac{562115611087}{46487926782} {S_3'}^2
    +\frac{762093077}{66364056} {S_2'} S'_4- \frac{1990577}{47223} S'_6 -\\-
    \frac{4217456129563}{116219816955}\cdot \frac{S'_3 S'_5}{S'_2}
    - \frac{14623247}{1301256} \cdot\frac{{S_4'}^2}{S'_2}+\frac{2359787}{31482}
    \frac{S'_8}{S'_2}.
\end{multline*}

The above root $S_6''$ must satisfy all other 13 equations, 
namely, the equations from 15 to 26 and 13th. 
Substituting this expression for $S_6''$ into 13 equations we obtain 13 polynomial
relations on $S_3'$, \dots, $S_{12}'$. They are cumbersome: the first one
is a sum of 10 monomials, the last one is a sum of 130 monomials.

The following observation simplifies calculations slightly.
It's easy to check that 13th equation linearly depends on $S_6$. 
If it has two different positive roots, then the coefficient of $S_6$ in it vanishes.
This allows us to express $S_7'$ in terms of $S_3'$, $S_4'$, $S_5'$ and
decrease the number of variables by 1. 
$$ 
S_7'=-\frac{1494661249487}{4501080325368} S_3'{S_2'}^2
+\frac{217002961}{417230286} S_2' S_5' +\frac{3678199}{2599908} S_3'S_4'.  
$$
Applying this substitution, we obtain a system of 12 equations.
The Gr\"obner basis of this huge system is relatively small.
Maple shows that the system has only two solutions:
the solution which gives us the examples presented in theorem 1
and a ``self-dual'' solution for which $S_6''=S_6'$.


\begin{thebibliography}{99}

\bibitem{fomin}
D. V. Fomin и O. T. Izhboldin, {\it Sets of multiple sums}. 
Proc. of St.Petersburg Math. Soc., {\bf 3},  244--259 (1994). %(in Russian)

\bibitem{boman}
J. Boman and L. Svante, {\it Examples of non-uniqueness for the combinatorial
Radon transform modulo the symmetric group}. 
Math. Scand., {\bf 78}, 207--212, (1996).

\bibitem{ewell}
J. A. Ewell, { \it On the determination of sets by sets of sums of fixed
order}.  Canad. J. Math., {\bf 20}, 596--611, (1968).

\bibitem{gordon}
A. S. Fraenkel, B. Gordon and E. G. Straus, { \it On the determination of sets
by the sets of sums of a certain order}. 
Pacific J. Math., {\bf 12},  187--196, (1962).

\bibitem{Guy1994}
R. K. Guy, {\it Unsolved problems in number theory}.
Springer, New York (1994).

\bibitem{MacMahon}
P. A. MacMahon {\it Combinatorial analysis}. Chelsea, New York (1960).

\bibitem{Moser1957}
L. Moser, {\it Problem E1248}. 
Amer. Math. Monthly, {\bf 64}, 507, (1957).

\bibitem{straus}
J. L. Selfridge and E. G. Straus, { \it On the determination of numbers by
their sums of a fixed order}. 
Pacific J. Math., {\bf 8}, 847--856, (1958).

\end{thebibliography}
\end{document}